\documentclass[a4paper,12pt]{article}

\usepackage{amsmath,amssymb,amsthm}
\usepackage{amsfonts}

\numberwithin{equation}{section}

\newtheorem{theorem}{Theorem}[section]
\newtheorem{lemma}[theorem]{Lemma}
\newtheorem{e-proposition}[theorem]{Proposition}
\newtheorem{corollary}[theorem]{Corollary}
\newtheorem{e-definition}[theorem]{Definition\rm}
\theoremstyle{definition}
\newtheorem{remark}[theorem]{\it Remark\/}

\def\ignore#1{}
\def\eq{\begin{equation}}
\def\en{\end{equation}}
\def\eqa{\begin{eqnarray}}
\def\ena{\end{eqnarray}}
\def\eqs{\begin{eqnarray*}}
\def\ens{\end{eqnarray*}}

\def\re{{\mathbb R}}
\def\non{\nonumber}
\def\s{\sigma}

\def\d{\delta}
\def\Ref#1{(\ref{#1})}%
\def\Ref{\eqref}
\def\Eq{\ =\ }
\def\Le{\ \le\ }

\def\h{\eta}
\def\t{\tau}
\def\th{\theta}
\def\g{\gamma}

\def\G{\Gamma}

\def\e{\varepsilon}
\def\f{\phi}
\def\a{\alpha}

\def\nn{{\mathcal N}}

\def\pr{{\mathbb P}}

\def\ex{{\mathbb E}}
\def\ui{^{(1)}}
\def\ut{^{(2)}}

\def\aa{{\mathcal A}}
\def\law{{\mathcal L}}
\def\L{\Lambda}
\def\r{\rho}
\def\D{\Delta}
\def\m{\mu}
\def\giv{\,|\,}
\def\Blb{\left\{}
\def\Brb{\right\}}

\def\bone{{\mathbf 1}}
\def\nin{\noindent}
\def\half{{\textstyle\frac12}}
\def\quarter{{\textstyle\frac14}}
\def\third{{\textstyle\frac13}}
\def\Bl{\Bigl(}
\def\Br{\Bigr)}

\def\Def{\ :=\ }
\def\nt{{\lfloor nt \rfloor}}
\def\nu{{\lfloor nu \rfloor}}
\def\nv{{\lfloor nv \rfloor}}
\def\sn{\sum_{i=1}^n}
\def\sjn{\sum_{j=1}^n}
\def\sln{\sum_{l=1}^n}
\def\skn{\sum_{k=1}^n}
\def\sijn{\sum_{i,j=1}^n}
\def\sjkn{\sum_{j,k=1}^n}

\def\sint{\sum_{i=1}^\nt}
\def\sjnu{\sum_{j=1}^\nu}
\def\un{^{(n)}}
\def\ta{\tilde a}
\def\tf{\tilde f}
\def\ts{\tilde s}
\def\tS{\widetilde S}

\def\tL{\widetilde \L}
\def\var{{\rm Var\,}}
\def\ba{{\bar a}}
\def\p{\pi}
\def\cov{{\rm Cov\,}}
\def\swd{\sup_{w\in D}}
\def\uk{^{(k)}}
\def\nw{\|w\|}
\def\ep{\hfill$\Box$\\[2ex]}
\def\bX{{\overline X}}
\def\Ge{\ \ge\ }
\def\intinti{\int_{I_n(t)\times I}}
\def\intiti{\int_{[0,t]\times I}}

\def\slogn{\sqrt{\log n}}
\def\te{{\tilde\e}}
\def\b{\beta}
\def\etal{{\it et al.\/}}
\def\Var{\var}

\newcommand\bigpar[1]{\bigl(#1\bigr)}
\newcommand\Bigpar[1]{\Bigl(#1\Bigr)}

\newcommand\qqw{^{-1/2}}

\newcommand\qw{^{-1}}

\newcommand\hY{\widehat Y}
\newcommand\hZ{\widehat Z}
\newcommand\hgs{\widehat \sigma}

\newcommand\bigset[1]{\ensuremath{\bigl\{#1\bigr\}}}

\begin{document}

\title{A functional combinatorial central limit theorem}

\author{A.\ D.\ Barbour\footnote{Angewandte Mathematik, Universit\"at Z\"urich, 
Winterthurertrasse 190, CH-8057 Z\"URICH; 
work supported in part by Schweizerischer Nationalfonds Projekte Nr.\ 
20--107935/1 and 20--117625/1.}\ \ and 
S.\ Janson\footnote{Matematiska institutionen, Uppsala universitet,
Box 480, SE--751~06 UPPSALA.\hfil\break
Research done while both authors visited the Institut Mittag--Leffler,
Djursholm, Sweden.}  
\\ 
Universit\"at Z\"urich and Uppsala Universitet}

\date{July 2, 2009}
\maketitle{}

\begin{abstract}
The paper establishes a functional version of the Hoeffding 
combinatorial central limit theorem.  First, a
pre-limiting Gaussian process approximation is defined, and is
shown to be at a distance of the order of the Lyapounov
ratio from the original random process.  Distance is
measured by comparison of expectations of smooth functionals
of the processes, and the argument is by way of Stein's method.  
The pre-limiting process is then shown,
under weak conditions, to converge to a Gaussian limit process.
The theorem is used to describe the shape of random
permutation tableaux.
\end{abstract} 


{\it Keywords:}\ {Gaussian process; combinatorial central limit theorem; 
\hfil\break\hglue2.58cm permutation tableau;
 Stein's method}


\setcounter{equation}{0}
\section{Introduction}
Let $a_0\un := (a_0\un(i,j),\,1\le i,j \le n)$, $n\ge1$, be a sequence of
real matrices.
Hoeffding's (1951) combinatorial central limit theorem asserts that
if~$\p$ is a uniform random permutation of $\{1,2,\ldots,n\}$,
then, under appropriate conditions, the distribution of the sum
$$
    S_0\un \Def \sn a_0\un(i,\p(i)) ,
$$
when centred and normalized, converges to the standard normal
distribution.   The centring is usually accomplished by
replacing $a_0\un(i,j)$ with 
\eqs
   \ta\un(i,j) &:=& a_0\un(i,j) - \ba_0\un(+,j) -  \ba_0\un(i,+)
      +  \ba_0\un(+,+),
\ens      
where
\eqs
    \ba_0\un(+,j) &:=& n^{-1}\sn a_0(i,j); \quad
    \ba_0\un(i,+) \Def n^{-1}\sjn a_0(i,j); \\
    \ba_0\un(+,+) &:=& n^{-2}\sn a_0(i,j).
\ens   
This gives $\tS\un = S_0\un - \ex S_0\un$, and the variance 
$\var \tS\un=\var S_0\un$
is then given by
$$
   \{\ts\un(a)\}^2 \Def (n-1)^{-1}\sijn  \{\ta\un(i,j)\}^2.
$$
Bolthausen~(1984) proved the analogous Berry--Esseen
theorem: that, for any $n \times n$ matrix~$a$, 
$$
    \sup_{x \in \re}|\pr[S_0 - m(a) \le x\ts(a)] - \Phi(x)| \Le C \tL(a),
$$
for a universal constant~$C$, where 
\begin{gather}
    S_0 := \sn a_0(i,\p(i)),
\qquad m(a) \Def n^{-1}\sijn a_0(i,j) \Eq \ex S_0,
\nonumber
\\\label{tsa}
     \ts^2(a) := (n-1)^{-1}\sijn \ta^2(i,j) \Eq \var S,   
  \end{gather}
(we tacitly assume $n\ge2$ when necessary) and 
$$
    \tL(a) \Def  \frac1{n\ts^{3}(a)}\sijn |\ta(i,j)|^3
$$
is the analogue of the Lyapounov ratio.    

In this paper, we begin by proving a functional version of
Bolthausen's theorem, again
with an error expressed in terms of a Lyapounov ratio. 
When centring the functional version $S_0(t) := \sint a_0(i,\p(i))$,
$0\le t\le 1$,
it is however no longer natural to make the double standardization that
is used to derive $\ta$ from~$a_0$.  Instead, we shall at each step
centre the random variables $a_0(i,\p(i))$ individually by
their means $\ba_0(i,+)$.  Equivalently, in what follows, we shall
work with matrices~$a$ satisfying $\ba(i,+)= 0$ for all~$i$,
but with no assumption as to the value of $\ba(+,j)$.
For example, if we have $a_0(i,j) = b(i) + c(j)$, then $\ta(i,j) = 0$
for all $i,j$, and hence 
$S_0=\ex S_0=n\ba_0(+,+)=n(\bar b+\bar c)$ is a.s.~constant.  However,
we are interested instead in 
\[
   S(t) \Def \sint \{a_0(i,\p(i)) - \ba_0(i,+)\},
\]	 
giving $S(t) = \sint \{c(\p(i)) - \bar c\}$, a non-trivial process
with a Brownian bridge as natural approximation.

We thus, throughout the paper, define the matrix~$a$ by
\begin{equation}
  \label{aa0}
    a(i,j)\Def a_0(i,j)-\ba_0(i,+),
\end{equation}
so that $\ba(i,+)=0$.
Correspondingly, we define
\begin{equation*}
    S(t)\Def \sint a(i,\pi(i)) \Eq S_0(t)-\ex S_0(t).
\end{equation*}
We then normalize by a suitable factor $s(a)>0$, and write
\begin{equation}
  \label{Y-def1}
     Y(t)\Def s(a)\qw S(t) = s(a)\qw\bigpar{S_0(t)-\ex S_0(t)};
\end{equation}
this can equivalently be expressed as 
\eq\label{Y-def2}
     Y \Def Y(\p) \Def \frac1{s(a)}\sn a(i,\p(i)) J_{i/n},
\en
where $J_u(t) := \bone_{[u,1]}(t)$.
In Theorem~\ref{Th1}, we approximate the random function $Y$
by the Gaussian process
\eq\label{Z-def}
    Z \Def \sn W_i J_{i/n},
\en
in which the jointly Gaussian random variables $(W_i, 1\le i\le n)$
have zero means and covariances given by
\begin{equation}
   \label{covariances} 
  \begin{aligned}
   \var W_i &\Eq \frac1{ns^2(a)}\sln a^2(i,l) \ =:\ \s_{ii}; \\
   \cov(W_i,W_j) &\Eq -\frac1{n(n-1)s^2(a)}\sln a(i,l)a(j,l) \ =:\
   \s_{ij}, 
\qquad i\ne j.
  \end{aligned}
\end{equation}
A simple calculation shows that
$\cov\bigpar{a(i,\pi(i)),a(j,\pi(j))}=s^2(a)\s_{ij}$ for all $i,j$, and
thus the covariance structures of the processes $Y$ and~$Z$ are
identical.
The error in the approximation
is expressed in terms of a probability metric defined 
in terms of comparison of expectations of certain smooth
functionals of the processes, and it is bounded by a multiple
of the Lyapounov ratio
\eq\label{L(a)-def}
   \L(a) \Def \frac1{ns^3(a)}\sijn |a(i,j)|^3.
\en   

The normalization factor $s(a)$ may be chosen in several ways. One 
obvious possibility is to choose $s(a)=\ts(a)$ defined in \eqref{tsa},
which makes $\var Y(1)=\var Z(1)=1$.
At other times this is inappropriate; for example, as seen above,
$\ts(a)$ may vanish, although we have a non-trivial Brownian bridge asymptotic.
A canonical choice of normalization is 
\eq\label{s(a)-def}
    s^2(a) \Def \frac1{n-1} \sijn a^2(i,j),     
\en 
or, for simplicity,  $n\qw\sijn a^2(i,j)$, which makes no difference
asymptotically.
In the special case where $\ba(+,j) = 0$ for each~$j$,
as with the matrix~$\ta$, this gives $s^2(a)=\ts^2(a)$, so 
$\Var Y(1)=\Var Z(1)=1$, 
but in general this does not hold.
In specific applications, some other choice may be more convenient. We
thus state our main results for an arbitrary normalization.

In most circumstances, such an approximation by $Z=Z\un$ depending on $n$ 
is in itself not
particularly useful; one would prefer to have some fixed, and
if possible well-known limiting approximation.  This requires
making additional assumptions about the sequence of 
matrices~$a\un$ as $n \to \infty$.  In extending
Bolthausen's theorem, it is enough to assume that $\tL\un(a)
\to 0$, since the approximation is already framed in terms
of the standard normal distribution.  For functional approximation,
even if we had standardized to make $\var Y(1) = 1$, we would
still have to make some further assumptions about the~$a\un$,
in order to arrive at a limit.  A natural choice would be to
take $a\un(i,j) := \a(i/n,j/n)$ for a continuous function
$\a\colon [0,1]^2 \to \re$ which does not depend on~$n$.
We shall make a somewhat weaker assumption, enough
to guarantee that the covariance function of~$Z\un$ converges
to a limit, which itself determines a limiting Gaussian process. 
The details are given in Theorem~\ref{Th2}.  Note that we
require that $\L\un(a) \log^2 n \to 0$ for process convergence,
a slightly stronger condition than might have been expected.
This is as a result of the method of proof, using the approach
in Barbour~(1990), in which the probability metric used for
approximation is not obviously strong enough to metrize weak
convergence in the Skorohod topology. Requiring the rate of
convergence of~$\L\un(a)$ to zero to be faster than
$1/\log^2 n$ is however enough to ensure that weak convergence also
takes place: see Proposition~\ref{prop1}.  

The motivation for proving the theorems comes from the study of
permutation tableaux.  In Section~\ref{example}, we show that the
boundary of a random permutation tableau, in the limit as its
size tends to infinity, has a particular shape, about which
the random fluctuations are approximately Gaussian.  The main
tool in proving this is Theorem~\ref{Th2}, applied to the
matrices $a_0\un(i,j) := \bone_{\{i\le j\}}$.

\setcounter{equation}{0}
\section{The pre-limiting approximation}
We wish to show that the distributions of the processes
$Y$ and~$Z$ of \Ref{Y-def2} and~\Ref{Z-def} are close.
To do so, we adopt the approach in Barbour~(1990). We
let~$M$ denote the space of all twice Fr\'echet differentiable
functionals $f\colon D := D[0,1] \to \re$ for which the norm
\eqa
  &&\|f\|_M \Def \swd\{|f(w)|/(1+\nw^3)\} +\swd\{\|Df(w)\|/(1+\nw^2)\} 
              \label{norm-def}\\   
   &&\mbox{}\quad + \swd\{\|D^2f(w)\|/(1+\nw)\} 
       + \sup_{w,h\in D}\{\|D^2f(w+h) - D^2f(w)\|/\|h\|\}
  \non
\ena
is finite; here, $\|\cdot\|$ denotes the supremum norm on~$D$, and
the norm of a (symmetric)  
$k$-linear form~$B$ on function in~$D$ is defined to
be $\|B\| := \sup_{h\in D\colon\,\|h\|=1} |B[h\uk]|$, where~$h\uk$
denotes the $k$-tuple $(h,h,\ldots,h)$.  Our aim is to show that
$|\ex g(Y) - \ex g(Z)|$ is small for all $g\in M$.  We do this by
Stein's method, observing that, for any $g\in M$, 
there exists a function $f \in M$ satisfying
\eq\label{Stein-eqn}
  g(w) - \ex g(Z) \Eq (\aa f)(w) \Def 
       -Df(w)[w] + \sijn \s_{ij} D^2f(w)[J_{i/n},J_{j/n}],
\en
and that 
\eq\label{C0-def}
      \|f\|_M \le C_0\|g\|_M,
\en 
where~$C_0$ does not depend on the choice of~$g$:
see, for example, Barbour (1990, (2.24), Remark 7 after Theorem 1 
and the remark following Lemma 3.1). 
Hence it is enough to prove that $|\ex(\aa f)(Y)| \le \e \|f\|_M$
for all $f\in M$ and for some small~$\e$.

\begin{theorem}\label{Th1}
Let $Y = Y(\p)$ and~$Z$ be defined as in \Ref{Y-def2} and~\Ref{Z-def},
with~$\p$ a uniform random permutation of $\{1,2,\ldots,n\}$, and
$\L(a)$ as in~\Ref{L(a)-def}, 
for some $n\times n$ matrix $a(i,j)$ with $\ba(i,+)=0$ and some $s(a)>0$.
Then there exists a universal constant~$K$ such that, for all $f \in M$,
$$
   | \ex(\aa f)(Y)| \Le K \L(a) \|f\|_M.
$$
Thus, for all $g\in M$,
$$
    |\ex g(Y) - \ex g(Z)| \Le C_0 K \L(a) \|g\|_M,
$$
with~$C_0$ as in~\Ref{C0-def}.
\end{theorem}

\proof
We begin by noting that
\eq\label{4}
   \ex Df(Y)[Y] \Eq \frac1{s(a)} \sn \ex\{X_i Df(Y)[J_{i/n}]\},
\en
where $X_i := a(i,\p(i))$.  We then write
\eq\label{5}
   \ex\{X_i Df(Y)[J_{i/n}]\} \Eq \frac1n \sln a(i,l) 
       \ex\{Df(Y(\p))[J_{i/n}] \giv \p(i) = l\}.
\en
Now realize~$\p'$ with the distribution $\law(\p\giv\p(i)=l)$ by
taking~$\p$ to be a uniform random permutation, and setting
\eqs
  &&\p' \Eq \p, \hspace{3.06in} \mbox{if}\quad \p(i) \Eq l;\\
  &&\p'(i) \Eq l;\quad \p'(\p^{-1}(l)) \Eq j;\quad \p'(k) \Eq \p(k),\ 
     k\notin \{i,\p^{-1}(l)\}, \\
	&&\hspace{3.7in}	  \mbox{if}\quad \p(i) \Eq j\ \neq\ l.
\ens
This gives
\eq\label{6}
  Y(\p') \Eq Y(\p) + \D_{il}(\p) \ =:\ Y'(\p),
\en
where
\eq\label{Delta-def}
  s(a)\D_{il}(\p) \Def \{a(i,l) - a(i,\p(i))\}J_{i/n}
   + \{a(\p^{-1}(l),\p(i)) - a(\p^{-1}(l),l)\}J_{\p^{-1}(l)/n},
\en
and $Y'(\p)$ has the distribution $\law(Y(\p) \giv \p(i)=l)$. 
Hence, putting \Ref{6} into~\Ref{5}, it follows that
\eq\label{8}
  \frac1{s(a)} \ex\{X_i Df(Y)[J_{i/n}]\} \Eq \frac1{ns(a)} \sln a(i,l) 
       \ex\{Df(Y(\p)+ \D_{il}(\p))[J_{i/n}]\}.
\en

Using Taylor's expansion, and recalling the definition~\Ref{norm-def}
of $\|\cdot\|_M$, we now have
\eqa
  \lefteqn{|\ex\{Df(Y+ \D_{il})[J_{i/n}]\} - \ex\{Df(Y)[J_{i/n}]\}
    - \ex\{D^2f(Y)[J_{i/n},\D_{il}]\}| } \non\\
   &&\Le \|f\|_M \ex\|\D_{il}\|^2, 
     \phantom{HHHHHHHHHHHHHHHHHHHHH} \label{9}
\ena
where, from \Ref{Delta-def},
\eq\label{9a}
  \|\D_{il}(\p)\| \Le \{s(a)\}^{-1}\{|a(i,l)| + |a(i,\p(i))| + |a(\p^{-1}(l),\p(i))|
     + |a(\p^{-1}(l),l)|\}.
\en
Laborious calculation now shows that
\eq\label{10}
  \frac1{ns(a)} \sn\sln |a(i,l)|\ex\|\D_{il}\|^2 \Le C_1\,\frac1{ns^3(a)}
     \sn\sln |a(i,l)|^3 \Eq C_1 \L(a),
\en
for a universal constant~$C_1$; for instance,
\eqs
  \lefteqn{\frac1{ns(a)} \sn\sln |a(i,l)| \frac1{s^2(a)} 
	        \ex|a(i,\p(i))a(\p^{-1}(l),\p(i))|} \\
	&&\Le \frac1{ns^3(a)}\sn\sln |a(i,l)|\Bigl\{ \frac1n a^2(i,l) +
	   \frac1{n(n-1)} \sum_{j\ne l}\sum_{k\ne i}|a(i,j)a(k,j)| \Bigr\} \\
	&&\Le \frac1{ns^3(a)}\sn\sln \Bigl\{\frac1n |a(i,l)|^3 +	 
	   \frac1{n(n-1)} \sum_{j\ne l}\sum_{k\ne i} \third\{|a(i,l)|^3 +	
		     |a(i,j)|^3 + |a(k,j)|^3 \} \Bigr\} \\
	&&\Eq  \frac1{ns^3(a)} \sn\sln |a(i,l)|^3.			 			
\ens
Thus, in view of~\Ref{8}, when evaluating the right hand side of~\Ref{4},
we have
\eq\label{11}
   \ex Df(Y)[Y] \Eq \frac1{ns(a)} \sn\sln a(i,l)\bigl(\ex \{Df(Y)[J_{i/n}]\} +
      \ex\{ D^2f(Y)[J_{i/n},\D_{il}]\}\bigr) + \h_1, 
\en
where $|\h_1| \le C_1\L(a)\|f\|_M$.  

Now, because $\ba(i,+) = 0$, the first term on the right hand side 
of~\Ref{11} is zero, so we have only the second to consider.
We begin by writing
\eq\label{12}
   D^2f(Y)[J_{i/n},\D_{il}] \Eq D^2f(Y)[J_{i/n},\ex\D_{il}] 
        + D^2f(Y)[J_{i/n},\D_{il} - \ex\D_{il}].
\en
{}From~\Ref{Delta-def}, it follows easily that
\eqa
   \ex\{D^2f(Y)[J_{i/n},\ex\D_{il}]\} 
       &=& \{s(a)\}^{-1} a(i,l) \ex\{D^2f(Y)[J_{i/n}\ut]\}  \label{13} \\
   &&\mbox{}
   - \frac1{(n-1)s(a)}\sum_{r\neq i} a(r,l) \ex\{D^2f(Y)[J_{i/n},J_{r/n}]\}.
   \non
\ena
Substituting this into~\Ref{11} gives a contribution to $\ex Df(Y)[Y]$ of
\eqa
   \f_1 &:=& \frac1{ns^2(a)}\sn\sln a^2(i,l) \ex\{D^2f(Y)[J_{i/n}\ut]\}\non \\
   &&\mbox{}\quad -  \frac1{n(n-1)s^2(a)}\sn\sln a(i,l) \sum_{r\neq i} a(r,l)
     \ex\{D^2f(Y)[J_{i/n},J_{r/n}]\} \non\\
   &=& \sn \s_{ii}\ex\{D^2f(Y)[J_{i/n}\ut]\} + \sn\sum_{r\neq i}\s_{ir}
     \ex\{D^2f(Y)[J_{i/n},J_{r/n}]\}, \label{14}
\ena 
from \Ref{covariances}, and thus, from \Ref{Stein-eqn}, \Ref{11} and~\Ref{12},
and noting that \Ref{14} cancels the second term in \Ref{Stein-eqn}, 
\eq\label{14a}
  | \ex(\aa f)(Y)| \Le |\h_1| + |\h_2|,
\en
where
\eq\label{14b}
   |\h_2| \Le  \frac1{ns(a)} \sn\sln |a(i,l)|\,
      |\ex\{ D^2f(Y)[J_{i/n},\D_{il} - \ex\D_{il}]\}|.
\en
It thus remains to find a bound for this last expression.

To address this last step, we write
\eqs
   \lefteqn{\ex\{ D^2f(Y)[J_{i/n},\D_{il} - \ex\D_{il}]\} }\\
  &&\Eq \sjkn p_{jk}\ex\{ D^2f(Y)[J_{i/n},\D_{il} - \ex\D_{il}] 
        \giv \p(i)=j,\p^{-1}(l)=k\},
\ens
where $p_{jk} := \pr[\p(i)=j,\p^{-1}(l)=k]$; note that $p_{li} = 1/n$,
and that $p_{jk} = 1/n(n-1)$ for $j\ne l$, $k\ne i$.  We then observe that, 
much as for~\Ref{6},
\eq\label{17}
   Y''(\p) \Def Y(\p) + \D'_{il;jk}(\p) \ \sim\ \law(Y(\p) \giv \p(i)=j,\p^{-1}(l)=k),
\en
where, for $j\ne l$, $k\ne i$,
\begin{equation}\label{18}
  \begin{split}
s(a)  \D'_{il;jk}(\p) &:= 
  \bigl\{[a(i,j) - a(i,\p(i))]J_{i/n} 
                     + [a(k,l) - a(k,\p(k))]J_{k/n} \\
     &\hskip3em + [a(\p^{-1}(l),\p(k)) - a(\p^{-1}(l),l)]J_{\p^{-1}(l)/n}\\
     &\hskip3em   + [a(\p^{-1}(j),\p(i)) - a(\p^{-1}(j),j)]J_{\p^{-1}(j)/n} \bigr\}
               \bone_{\{\p(i)\ne l,\,\p(k)\ne j\}} \\
      &\quad + 
\bigl\{[a(i,j) - a(i,\p(i))]J_{i/n} 
                     + [a(k,l) - a(k,j)]J_{k/n}\\ 
     &\hskip3em  + [a(\p^{-1}(l),\p(i)) -
                     a(\p^{-1}(l),l)]J_{\p^{-1}(l)/n} \bigr\}
        \bone_{\{\p(i)\ne l,\,\p(k)= j\}} \\
      &\quad + 
\bigl\{[a(i,j) - a(i,l)]J_{i/n} 
                     + [a(k,l) - a(k,\p(k))]J_{k/n}\\ 
     &\hskip3em  + [a(\p^{-1}(j),\p(k)) - a(\p^{-1}(j),j)]J_{\p^{-1}(j)/n} \bigr\}
        \bone_{\{\p(i)= l\}},  
  \end{split}
\end{equation}
and
\eq
  s(a)\D'_{il;li}(\p) \Def [a(i,l) - a(i,\p(i))]J_{i/n}
    + [a(\p^{-1}(l),\p(i)) - a(\p^{-1}(l),l)]J_{\p^{-1}(l)/n}. \label{18a}
\en
Then $\D_{il} = \D_{il}(\p(i),\p^{-1}(l))$ is measurable with respect to 
$\s(\p(i),\p^{-1}(l))$, and 
\eqa
  \lefteqn{\sjn\skn p_{jk} \ex\{D^2f(Y)[J_{i/n},\D_{il} - \ex\D_{il}] \giv 
             \p(i)=j,\p^{-1}(l)=k\} }\non\\
  &&\Eq \sjn\skn p_{jk} \ex\{D^2f(Y+\D'_{il;jk})
          [J_{i/n},\D_{il}(j,k) - \ex\D_{il}]\}    \non\\
  &&\Eq   \sjn\skn p_{jk} \ex\{D^2f(Y)[J_{i/n},\D_{il}(j,k) - \ex\D_{il}] \} 
                    \non\\
  &&\mbox{}\qquad
    + \sjn\skn p_{jk}\ex\{D^2f(Y+\D'_{il;jk})
                   [J_{i/n},\D_{il}(j,k) - \ex\D_{il}] \non\\
  &&\mbox{}\hspace{3.5cm}  - D^2f(Y)[J_{i/n},\D_{il}(j,k) - \ex\D_{il}]\}. 
         \label{19}
\ena
Now, since $\sjn\skn p_{jk}\D_{il}(j,k) = \ex\D_{il}$, the first term
in~\Ref{19} is zero, by bilinearity.  For the remainder, we have
\eqa
  \lefteqn{\| D^2f(Y+\D'_{il;jk})[J_{i/n},\D_{il}(j,k) - \ex\D_{il}] 
          - D^2f(Y)[J_{i/n},\D_{il}(j,k) - \ex\D_{il}] \| } \non\\
   &&\Le \|f\|_M \|\D'_{il;jk}\|\{\|\D_{il}(j,k)\| + \|\ex\D_{il}\|\},
   \phantom{HHHHHHHHHHHH}
   \label{20}       
\ena
so that, from~\Ref{14b}, 
\eq\label{21}
  |\h_2|
     \Le \|f\|_M  \frac1{ns(a)} \sn\sln |a(i,l)|\,
   \sjn\skn p_{jk} \ex \|\D'_{il;jk}\| \{\|\D_{il}(j,k)\| + \|\ex\D_{il}\|\}.
\en  
Here, from \Ref{Delta-def}, \Ref{9a} and~\Ref{18}, each of the norms
can be expressed as $1/s(a)$ times a sum of elements of~$|a|$. 
Another laborious calculation shows that indeed
$$
   |\h_2| \Le C_2 \L(a) \|f\|_M,
$$
and the theorem is proved.
\ep

\setcounter{equation}{0}
\section{A functional limit theorem}
The pre-limiting approximation is simpler than the original
process, inasmuch as it involves only jointly Gaussian random
variables with prescribed covariances.  However, if the 
matrix~$a$ can be naturally imbedded into a sequence~$a\un$
exhibiting some regularity as~$n$ varies, and if~$n$ is large,
it may be advantageous to look for an $n$-independent
limiting approximation, in the usual sense of weak convergence.
Unfortunately, the approximation given in Theorem~\ref{Th1}
is not naturally compatible with weak convergence with respect
to the Skorohod metric, and something extra is needed.  With
this in mind, we prove the following extension of
Theorem~2 of Barbour~(1990).  To do so, we introduce the
class of functionals $g \in M^0 \subset M$ for which
$$
  \|g\|_{M^0} \Def \|g\|_M + \swd|g(w)| + \swd\|D g(w)\| + \swd\|D^2g(w)\| 
   \ <\ \infty.
$$
   
\begin{e-proposition}\label{prop1}
Suppose that, for each $n\ge1$, the random element~$Y_n$
of~$D := D[0,1]$ 
is piecewise constant, with intervals of constancy of length 
at least~$r_n$. Let $Z_n$, $n\ge 1$, be random elements 
of~$D$ 
converging weakly in~$D$ to a random 
element~$Z$ of $C[0,1]$.  Then, if
\eq\label{g-bnd}
   |\ex g(Y_n) - \ex g(Z_n)| \Le C\t_n \|g\|_{M^0}
\en
for each $g\in M^0$, and if $\t_n \log^2(1/r_n) \to 0$
as $n\to\infty$, then $Y_n \to Z$ in~$D$.
\end{e-proposition}

\proof
First note that, by Skorohod's representation theorem, 
we may assume that the
processes $Z_n$ and~$Z$ are all defined on the same probability
space, in such a way that $Z_n \to Z$ in~$D$ a.s.\ as
$n \to \infty$.  Since~$Z$
is continuous, this implies that $\|Z_n-Z\| \to 0$ a.s.

As in the proof of Barbour~(1990, Theorem~2), it is enough
to show that 
\eq\label{Y-to-Z}
   \pr[Y_n \in B] \to \pr[Z \in B]
\en
for all
sets~$B$ of the form $\bigcap_{1\le l\le L}B_l$, where
\hbox{$B_l = \{w\in D\colon \|w-s_l\| < \g_l\}$} for $s_l\in C[0,1]$, 
and $\pr[Z \in \partial B_l] = 0$.  To do so, we approximate the
indicators $I[Y_n \in B_l]$ from above and below by
functions from a family $g := g\{\e,p,\r,\h,s\}$ 
in~$M^0$, and use~\Ref{g-bnd}. We define
$$
  g\{\e,p,\r,\h,s\}(w) \Def \f_{\r,\h}(h_{\e,p}(w-s)),
$$
where
$$
   h_{\e,p}(y) \Def \Bl \int_0^1(\e^2 + y^2(t))^{p/2}\,dt \Br^{1/p}
      \ =:\ \|(\e^2 + y^2)^{1/2}\|_p,
$$
and $\f_{\r,\h}(x) := \f((x-\r)/\h)$, for $\f\colon \re^+ \to [0,1]$
non-increasing, three times continuously differentiable, and
such that $\f(x) = 1$ for $x\le 0$ and $\f(x) = 0$ for $x \ge 1$.
Note that each such function~$g$ is in~$M^0$, and that
$\|g\|_{M^0} \le C'p^2\e^{-2}\h^{-3}$ for a constant~$C'$ not
depending on $\e,p,\r,\h,s$, and that the same is true for
finite products of such functions, if the largest of the $p$'s
and the smallest of the $\e$'s and $\h$'s is used in the norm 
bound.    
 
Now, if $x \in B_l$, it follows that $g_l(x) = 1$, for
$$
   g_l \Def g\{\e\g_l,p,\g_l(1+\e^2)^{1/2},\h,s_l\},
$$
for all $\e,p,\h$.   Hence, for all $\e,p,\h$,
\eq\label{Y-below-1}
  \pr\Bigl[Y_n \in \bigcap_{1\le l\le L}B_l\Bigr]
   \Le \ex\Bigl\{\prod_{i=1}^L g_l(Y_n) \Bigr\}   
   \Le \ex\Bigl\{\prod_{i=1}^L g_l(Z_n) \Bigr\}
     + C \t_n\,C'_B p^2(\e\g)^{-2}\h^{-3},
\en
where $\g := \min_{1\le l\le L}\g_l$. Then, by
Minkowski's inequality,
$$
   h_{\e,p}(Z-s_l) \Le h_{\e,p}(Z_n-s_l) + \|Z_n-Z\|_p
                   \Le h_{\e,p}(Z_n-s_l) + \|Z_n-Z\|.
$$
Hence, if $p_n \to \infty$ as $n \to \infty$ and $\e$ is
fixed, 
$$
  \liminf_{n \to \infty}h_{\e,p_n}(Z_n-s_l)                   
   \ \ge\ \liminf_{n \to \infty}\{h_{\e,p_n}(Z-s_l) - \|Z_n-Z\|\}
   \Eq \|(\e^2 + |Z-s_l|^2)^{1/2}\| \ a.s.
$$
It thus follows that, if $\|Z-s_l\| > \g_l$, and if $\h_n \to 0$
as $n \to \infty$, then
$$
    \liminf_{n \to \infty}\{h_{\e\g_l,p_n}(Z_n-s_l) - \h_n\}                   
   \ \ge\ \|(\e^2\g_l^2 + |Z-s_l|^2)^{1/2}\| \ >\ \g_l(1+\e^2)^{1/2} 
       \ a.s.,
$$
and so $g_{ln}(Z_n) = 0$ for all~$n$ sufficiently large,
where 
\[
   g_{ln} := g\{\e\g_l,p_n,\g_l(1+\e^2)^{1/2},\h_n,s_l\}. 
\]
Applying Fatou's lemma to $1 - \prod_{l=1}^L g_{ln}(Z_n)$, 
and because $\pr[Z \in \partial B_l] = 0$ for each~$l$, 
we then have,
\eqs
  \limsup_{n \to \infty} \ex\Bigl\{\prod_{i=1}^L g_{ln}(Z_n) \Bigr\}
  &\le& \ex\Bigl\{\limsup_{n \to \infty}\prod_{i=1}^L g_{ln}(Z_n) \Bigr\}\\      
  &\le& \ex\Bl \prod_{i=1}^L\bone\{\|Z-s_l\| \le \g_l\}\Br \Eq \pr[Z \in B].
\ens  
Thus, letting $p_n \to\infty$ and $\h_n \to 0$ in such a way
that $\t_n p_n^2\h_n^{-3} \to 0$, it follows from~\Ref{Y-below-1}
that $\limsup_{n \to \infty}\pr[Y_n \in B] \le \pr[Z \in B]$,
and we have proved one direction of~\Ref{Y-to-Z}.

For the other direction, fix $\th > 0$ small, and 
let~$\d>0$ be such that, if $\|Y_n-s_l\| \ge \g_l$, then 
\begin{equation}
  \label{geom}
  {\rm leb}\bigl\{t\colon\,|Y_n(t)-s_l| \ge \g_l(1-\th)\bigr\} \ge 
	   \bigl(\d\wedge\half r_n\bigr).
\end{equation}
Such a~$\d$ exists, because the collection $(s_l,\,1\le l\le L)$ is
uniformly equicontinuous, and because the functions~$Y_n$ are piecewise
constant on intervals of length at least~$r_n$.  Hence, for such~$Y_n$,
\[
   h_{\e\g_l,p}(Y_n-s_l) 
      \Ge \g_l\{\e^2+(1-\th)^2\}^{1/2}\bigl(\d\wedge\half r_n\bigr)^{1/p},
\]
and thus $g^*_l(Y_n) = 0$, where, for any $p$ and~$\h$,
\[
   g^*_l \Def g\bigl\{\e\g_l,p,
     \g_l(\e^2+(1-\th)^2)^{1/2}\bigl(\d\wedge\half r_n\bigr)^{1/p}-\h,\h,s_l\bigr\}.
\]
Thus, for any $p$ and~$h$, $I[Y_n \in B_l] \ge g_l^*(Y_n)$, and hence
\eq\label{Y-above-1}
   \pr\Bigl[Y_n \in \bigcap_{1\le l\le L}B_l\Bigr]
   \Ge \ex\Bigl\{\prod_{i=1}^L g_l^*(Y_n) \Bigr\}   
   \Ge \ex\Bigl\{\prod_{i=1}^L g_l^*(Z_n) \Bigr\}
     - C \t_n\,C'_B p^2(\e\g)^{-2}\h^{-3}.
\en

Now suppose that
$\|Z-s_l\| < \g_l(1-\th)$.  Then there exists an~$\a > 0$ such
that a.s.\ $\|Z_n-s_l\| < \g_l(1-\th) - \a$ for all~$n$ sufficiently
large. This in turn implies that
\eqs
  h_{\e\g_l,p_n}(Z_n-s_l) &\le& \{\e^2\g_l^2 + \|Z_n-s_l\|^2\}^{1/2}
    \Le \g_l\{\e^2 + (1-\th-\a\g_l^{-1})^2\}^{1/2} \\
   &<&  \g_l\{\e^2 + (1-\th)^2\}^{1/2}\bigl(\d\wedge\half r_n\bigr)^{1/p_n} - \h_n
\ens
for all~$n$ large enough, if $\h_n\to0$ and $p_n\to\infty$ in such a
way that $r_n^{1/p_n} \to 1$.  This in turn implies that $g_{ln}^*(Z_n)
= 1$ for all~$n$ large enough, where
\eq\label{gln-star-def}
   g_{ln}^* \Def g\bigl\{\e\g_l,p_n,
     \g_l(\e^2+(1-\th)^2)^{1/2}\bigl(\d\wedge\half r_n\bigr)^{1/p_n}
		       -\h_n,\h_n,s_l\bigr\}.
\en
Hence
\eq
   \ex\Bigl\{ \liminf_{n\to\infty}\prod_{i=1}^L g_{ln}^*(Z_n) \Bigr\}
   \Ge \pr\Bigl[ \bigcap_{1\le l\le L}\bigl\{\|Z-s_l\| < \g_l(1-\th)\bigr\} 
	 \Bigr].    \label{Y-above-3}
\en
Applying Fatou's lemma, and recalling~\Ref{Y-above-1}, we now have a.s.\
\eq
  \liminf_{n\to\infty}\pr\Bigl[Y_n \in \bigcap_{1\le l\le L}B_l\Bigr]  
   \Ge  \liminf_{n\to\infty}\ex\Bigl\{\prod_{i=1}^L g_{ln}^*(Z_n) \Bigr\}   
   \Ge \ex\Bigl\{ \liminf_{n\to\infty}\prod_{i=1}^L g_{ln}^*(Z_n) \Bigr\},
      \label{Y-above-2}
\en
provided that also $\t_n p_n^2 \h_n^{-3} \to 0$: this can be arranged by
judicious choice of $p_n \to \infty$ and $\h_n \to 0$ if, as assumed,
$\t_n \log^2(1/r_n) \to 0$.  
Hence, since $\th$ was chosen arbitrarily,
it follows from \Ref{Y-above-3} and~\Ref{Y-above-2} that
$$
  \liminf_{n\to\infty}\pr[Y_n \in B] \ge \pr[Z \in B],
$$
and the theorem is proved.  
\ep

\nin Note that, in Barbour~(1990, Theorem~2), restricting to functions~$g$
satisfying~(2.32) of that paper is not permissible:  the bound~\Ref{g-bnd}
is needed for functions in~$M_0$ that do not necessarily satisfy~(2.32).

\begin{remark} 
The assumption that $Y_n$ is piecewise constant can be relaxed to 
$Y_n$ being piecewise linear, with intervals of linearity of length
at least $r_n$; in particular, this allows processes~$Y_n$ obtained
by linear interpolation.
The only difference in the proof is that, if 
$\|Y_n-s_l\|\ge\g_l$, then $|Y_n(t_0)-s_l(t_0)|>(1-\th/4)\g_l$ for
some $t_0$. Thus, by the assumption on~$Y_n$ and the continuity of~$s_l$,
there exists an interval~$I_0$ of length at least
$l_n:=\frac12r_n\wedge\delta$, with $t_0$ as an endpoint,
on which~$Y_n$ is linear and $|s_l(t)-s_l(t_0)|<\th\g_l/4$.
A simple geometrical argument now shows that
$|Y_n(t)-s_l(t_0)|>(1-\th/2)\g_l$ in a subinterval of length at least
$\th l_n/8$, at one or other end of~$I_0$.
Hence, \Ref{geom} can be replaced by 
\begin{equation*}
  {\rm leb}\bigl\{t\colon\,|Y_n(t)-s_l| \ge \g_l(1-\th)\bigr\} \ge 
	   \tfrac\th{16}\bigl(\d\wedge r_n\bigr),
\end{equation*}
and the rest of the proof is the same.
\end{remark}

We now turn to proving a functional limit theorem for the sums derived
from  a sequence of matrices~$a\un$, $n\ge1$.  Supposing that $s\un(a) > 0$,
we define functions
\begin{equation}   \label{fun-def}
  \begin{aligned}
  f_n(t) &:= \frac1{n(s\un(a))^2}\sint \sln (a\un(i,l))^2;\\
  g_n(t,u) &:= \frac1{(ns\un(a))^2}\sint \sjnu\sln a\un(i,l)a\un(j,l),	
  \end{aligned}
\end{equation}
for $0\le t,u\le 1$.
Note that if we choose $s\un(a)$ by \Ref{s(a)-def}, then 
$f_n(1)=(n-1)/n\to1$. Conversely, if $f_n(1)$ converges to a 
limit $c>0$,
then $s\un(a)$ differs from the value in \Ref{s(a)-def} only by a
factor $c\qqw+o(1)$.

\begin{theorem}\label{Th2}
Suppose that $f_n \to f$ and $g_n \to g$ pointwise, with~$f$ continuous, 
and that $\L\un(a) \log^2 n \to 0$.  Then there exists a zero mean
continuous Gaussian 
process~$Z$ on $[0,1]$ with covariance function given by
\eq\label{b1c}
   \cov(Z(t),Z(u)) \Eq \s(t,u) \Def f(t\wedge u) - g(t,u),
\en
and $Y_n \to Z$ in $D[0,1]$.
\end{theorem}

\proof
Fix $n\ge2$.
We begin by realizing the random variables~$W_i\un$ as functions of
a collection $(X_{il},\,i,l\ge1)$ of independent standard normal
random variables.  Writing $\bX_l := n^{-1}\sn X_{il}$, we set
\eq\label{w-def}
   W_{il}\un \Def \frac1{s\un(a)\sqrt{n-1}}\, a\un(i,l)(X_{il} - \bX_l);\qquad
   W_i\un \Def \sln W_{il}\un.
\en
Direct calculation shows that, with $\d_{ij}$ the Kronecker delta,
\eqa
   \cov(W_i\un,W_j\un) &=& \sln \cov(W_{il}\un,W_{jl}\un) \label{cov-calc}\\ 
    &=& \sln \frac1{(n-1)(s\un(a))^2} a\un(i,l)a\un(j,l)(\d_{ij} - n^{-1}),
		\non
\ena
in accordance with~\Ref{covariances}, so we can set
\eq\label{Zn-def}
    Z_n \Def \sn W_i\un J_{i/n}.
\en
Now Theorem~\ref{Th1} shows that $|\ex\{g(Y_n) - g(Z_n)\}| \le
C\L\un(a) \|g\|_{M^0}$ for any $g\in M^0$; furthermore, the process~$Y_n$ is
piecewise constant on intervals of lengths~$1/n$, and, by assumption,
$\L\un(a) \log^2 n \to 0$.  Hence, in order to apply Proposition~\ref{prop1},
it is enough to show that $Z_n \to Z$ for a continuous Gaussian process.

Write $Z_n = Z_n\ui - Z_n\ut$, where 
\eq\label{Zn-split}
\begin{aligned}
     Z_n\ui(t) &\Def \frac1{s\un(a)\sqrt{n-1}} \sint \sln a\un(i,l)X_{il},
\\
     Z_n\ut(t) &\Def \frac1{s\un(a)\sqrt{n-1}} \sint \sln a\un(i,l)\bX_l.  
\end{aligned}
\en
The process~$Z_n\ui$ is a Gaussian process with independent increments,
and can be realized as $W(\tf_n(\cdot))$, where~$W$ is a standard Brownian
motion and $\tf_n(t) := nf_n(t)/(n-1)$.  Now~$f$ is continuous, by assumption, 
and each~$\tf_n$ is non-decreasing, so $\tf_n \to f$ uniformly on $[0,1]$, and 
hence $W(\tf_n(\cdot)) \to W(f(\cdot))$ in $D[0,1]$.  Since the latter process
is continuous, it follows that the sequence~$Z_n\ui$ is $C$-tight in $D[0,1]$. 

To show that $Z_n\ut$ is also $C$-tight, we use criteria from Billingsley~(1968).
For $0 \le t \le u \le 1$, it follows from~\Ref{Zn-split} and H\"older's
inequality that
\eqs
  \ex|Z_n\ut(u)-Z_n\ut(t)|^2 &=& \frac1{(n-1)(s\un(a))^2}\sln
      \Bl \sum_{i = \nt+1}^{\nu} a\un(i,l) \Br^2\,\frac1n \\
   &\le&  \frac1{n(n-1)(s\un(a))^2}(\nu-\nt)\sln \sum_{i = \nt+1}^{\nu} (a\un(i,l))^2\\
   &\le&  f_n(1)\frac{\nu-\nt}{n-1}.
\ens
Hence, since $Z_n\ut$ is Gaussian, we have
\eq\label{b7a}
  \ex|Z_n\ut(u)-Z_n\ut(t)|^4 \Eq 3(\ex|Z_n\ut(u)-Z_n\ut(t)|^2)^2 \Le 
     3\Bl f_n(1)\frac{\nu-\nt}{n-1} \Br^2.
\en
Thus, if $0 \le t \le v \le u \le 1$ and $u-t\ge1/n$, it follows that
\eqa
  &&\ex\Blb |Z_n\ut(v)-Z_n\ut(t)|^2 |Z_n\ut(u)-Z_n\ut(v)|^2 \Brb \non\\
   &&\qquad\Le \sqrt{\ex|Z_n\ut(v)-Z_n\ut(t)|^4\,\ex|Z_n\ut(u)-Z_n\ut(v)|^4} \non\\
  &&\qquad \Le 3 f^2_n(1)\Bl\frac{\nv-\nt}{n-1}\,\frac{\nu-\nv}{n-1} \Br 
       \Le 12 f_n^2(1)(u-t)^2; \phantom{HHH}\label{b7b}
\ena 
the inequality is immediate for $u-t < 1/n$, since then $\nv \in \{\nt,\nu\}$.

Now, for any $0\le t\le u\le 1$, we have
\[
   \cov(Z_n\ut(t),Z_n\ut(u)) \Eq \frac n{n-1}g_n(t,u) \ \to\ g(t,u).
\]
Hence there exists a zero mean Gaussian process~$Z\ut$ with covariance
function~$g$, and the finite dimensional distributions of~$Z_n\ut$
converge to those of~$Z\ut$.  By~\Ref{b7a} and Fatou's lemma,
$\ex|Z\ut(u)-Z\ut(t)|^4 \le 3f_n^2(1)(u-t)^2$ for any $0\le t\le u\le 1$,
so that, from Billingsley~(1968, Theorem~12.4), we may assume that
$Z\ut \in C[0,1]$.  From~\Ref{b7b} and Billingsley~(1968, Theorem~15.6),
it now follows that $Z_n\ut \to Z\ut$ in $D[0,1]$.  Thus $Z_n\ut$ is
$C$-tight also.

Now, since both $\{Z_n\ui\}$ and~$\{Z_n\ut\}$ are $C$-tight, so is their
difference $\{Z_n\}$.  From \Ref{fun-def} and~\Ref{cov-calc}, 
for $t,u\in[0,1]$,
\[
    \cov(Z_n(t),Z_n(u)) 
= \frac{n}{n-1} f_n(t\wedge u) - \frac{n}{n-1}g_n(t,u)
\ \to\ f(t\wedge u) - g(t,u),
\]
so that the finite dimensional distributions of~$Z_n$ converge to those
of a random element~$Z$ of $C[0,1]$ with covariance function~$\s(t,u)$,
as required.
\ep

\setcounter{equation}{0}
\section{Rate of convergence}
Under more stringent assumptions, the approximation of~$Z_n$ by~$Z$ can
be made sharper.  To start with, note that it follows from the
representation \Ref{w-def} and~\Ref{Zn-def} that~$Z_n$ can be written
as a two dimensional stochastic integral
\eq\label{Zn-rep}
   Z_n(t) \Eq \frac n{s\un(a)\sqrt{n-1}}\,\intinti \a_n(v,w)\,K(dv,dw)
\en
with respect to a Kiefer process~$K$,
where $I_n(t) := [0,n^{-1}\lfloor nt \rfloor]$, $I := [0,1]$ and
$\a_n(v,w) := a\un(\lceil nv\rceil,\lceil nw\rceil)$.
Recall that~the Kiefer process $K$ has covariance function
$\cov(K(v_1,w_1),K(v_2,w_2))=(v_1\wedge v_2-v_1v_2)(w_1\wedge w_2)$
and can be represented in the form $K(v,w) = W(v,w) - vW(1,w)$, where~$W$
is the two-dimensional Brownian sheet (Shorack \&~Wellner 1986, (5) p.~30
and Exercise~12, p.~32). Thus~$K$ is like a Brownian
bridge in~$v$, and a Brownian motion in~$w$.

In this section, we let $s(a\un)$ be given by \Ref{s(a)-def}. 
Hence if, for example, the functions~$\a_n$ converge in~$L_2$ to
a square integrable limit~$\a$ (not a.e.\ 0), then, 
\[
   \frac{n-1}{n^2}\{s\un(a)\}^2 \Eq \|\a_n\|_2^2\ \to\ \s^2_a 
     \Def \int_0^1\!dv \!\int_0^1\!dw \,\a^2(v,w) \Eq \|\a\|_2^2,
\]
and the limiting process~$Z$ can be represented as
\eq\label{Z-def-b}
   Z(t) \Eq \s_a^{-1}\intiti \a(v,w)\,K(dv,dw),
\en
enabling a direct comparison between $Z_n$ and~$Z$ to be made.
Since $\a_n \to_{L_2} \a$, it follows that 
\begin{equation} \label{fun-def-a}
  \begin{aligned}
   f_n(t) &\to\, f(t) \Def \s_a^{-2} \int_0^t\! dv \!\int_0^1\! dw \,\a^2(v,w);\\
	 g_n(t,u) &\to \,g(t,u) \Def \s_a^{-2} \int_0^t\! dv\! \int_0^u \!dx \!\int_0^1\! dw 
	    \,\a(v,w)\a(x,w),
  \end{aligned}
\end{equation}
with $f$ continuous, as required for Theorem~\ref{Th2}, and that~$Z$
has covariance function $\s(t,u)$ as defined in~\Ref{b1c}.
For the following lemma, we work under silghtly stronger assumptions.

\begin{lemma}\label{a-convergent-1}
Suppose that $\a_n \to \a$ in~$L_2$, where 
$\a$ is bounded and not a.e.~$0$,
and that, for some $0 < \b \le 2$, 
\eq\label{Fernique}
  |g(t,t) + g(u,u) - 2g(t,u)| \Le C_g^2|u-t|^\b,\qquad 0 \le t\le u\le 1.
\en	
Define 
$\a^+ := \|\a\|_\infty/\|\a\|_2 < \infty$ and 
$\e_n(v,w) := \|\a\|_2^{-1}\{\a_n(v,w) - \a(v,w)\}$.  
Then, for any $r > 0$, there is a constant $c(r)$ such that
\[
    \pr\bigl[\sup_{t\in I} |Z_n(t) - Z(t)| 
        > c(r)\bigl\{\|\e_n\|_2 + (\a^+  + C_g)n^{-(\b\wedge1)/2}\bigr\}\slogn\bigr] \Le n^{-r},
\]
where $Z$ is as defined in~\Ref{Z-def-b}.
\end{lemma}

\proof
Define $\te_n(v,w) := \frac n{s\un(a)\sqrt{n-1}}\a_n(v,w) - \s_a^{-1}\a(v,w)$.
We start by considering $t$ of the form $i/n$, $1\le i\le n$, so that
\[
   Z_n(t) - Z(t) \Eq \intiti \te_n(v,w)\,K(dv,dw).
\] 
{}From this and the representation $K(v,w) = W(v,w) - vW(1,w)$,
it follows that
$\max_{t\in I}\ex\{Z_n(t)-Z(t)\}^2 \le \|\te_n\|_2^2$, and hence, from 
the Borell--TIS maximal
inequality for Gaussian processes (Adler and Taylor~2007, Theorem~2.1.1), we have
\[
   \pr\Bigl[\max_{t\in n^{-1}\{1,2,\ldots,n\}}\Bigl|\intiti \te_n(v,w)\,K(dv,dw)\Bigr|
      > c_1(r)\|\te_n\|_2\slogn \Bigr] \Le \half n^{-r},
\]
if $c_1(r)$ is chosen large enough.  However, 
\[
   \te_n \Eq 
\frac{\a_n}{\|\a_n\|_2} - \frac{\a}{\|\a\|_2}, 
\]
from which it follows that
\[
  \|\te_n\|_2 \Le  2\|\e_n\|_2.
\]

It thus remains to consider the differences $Z_n(t) - Z(t)$ for~$t$ not of
the form $i/n$.  Between $n^{-1}\lfloor nt \rfloor$  and~$t$, the process~$Z_n$
remains constant, whereas~$Z$ changes; hence it is enough to control the
maximal fluctuation of~$Z$ over intervals of the form $[(i-1)/n,i/n]$,
$1\le i\le n$.  Here, we use  the Fernique--Marcus maximal inequality for 
Gaussian processes (Leadbetter~\etal\ 1983, Lemma~12.2.1), 
together with the inequality
\[
   |\s(u,u) + \s(t,t) - 2\s(t,u)| \Le C_g^2|t-u|^\b + (\a^+)^2|t-u|,
\]
to give the bound
\[
   \pr\Bigl[\max_{1\le i\le n}\sup_{(i-1)/n\le v\le i/n} 
     |Z(v) - Z((i-1)/n)| >  c_2(r)(C_g + \a^+) n^{-(\b\wedge1)/2}\slogn \Bigr]
         \Le \half n^{-r},
\]
if $c_2(r)$ is chosen large enough, and
the proof is now complete.
\ep

Note that, under the conditions of Lemma~\ref{a-convergent-1}, the
requirements  
for Theorem~\ref{Th2} are fulfilled, 
provided that $\L\un(a) \to 0$ fast enough.  This is true if also, for
instance,  
for some $c < \infty$,
$\|\a_n\|_\infty \le c\|\a\|_\infty$ for all~$n$, since then $\L\un(a)
\Le 2c\a^+n^{-1/2}$ 
for all~$n$ large enough.
Combining Theorems \ref{Th1} and~\ref{Th2}
with Lemma~\ref{a-convergent-1} then easily gives the following conclusions.

\begin{theorem}\label{a-convergent}
Under the conditions of Lemma~\ref{a-convergent-1}, and if also 
$\|\a_n\|_\infty/\|\a\|_\infty$ is bounded, then
$Y_n \to_d Z$ in $D[0,1]$, for~$Z$ as defined in~\Ref{Z-def-b},
and, for any functional $g\in M_0$, 
\eq\label{g-bnd-a}
   |\ex g(Y_n) - \ex g(Z)| \Le C\bigl\{\L\un(a) + n^{-1} +
	         \{\|\e_n\|_2 + (\a^+ + C_g) n^{-(\b\wedge1)/2}\}\slogn\bigr\}
	    \,\|g\|_{M^0},
\en
for some constant~$C$.
\end{theorem}

\proof
We note that
\[
    |\ex g(Y_n) - \ex g(Z)| \Le  |\ex g(Y_n) - \ex g(Z_n)| + \ex |g(Z_n) -  g(Z)|.
\]
The first term is bounded using Theorem~\ref{Th1}, whereas, for any $a>0$,
\eqs
    \ex |g(Z_n) -  g(Z)| &\le& 2\sup_{w\in D}|g(w)| \pr[\|Z_n-Z\|_\infty > a]
       + a\sup_{w\in D}\|Dg(w)\| \\
    &\le& \|g\|_{M_0} \{2\pr[\|Z_n-Z\|_\infty > a] + a\},
\ens
and the theorem follows by taking
$a = c(1)\bigl\{\|\e_n\|_2 + (\a^+  + C_g)n^{-(\b\wedge1)/2}\bigr\}\slogn$
and applying Lemma~\ref{a-convergent-1} with $r=1$.
\ep

\setcounter{equation}{0}
\section{The shape of permutation tableaux}\label{example}
We begin by studying the number of \emph{weak exceedances} in a uniform
random permutation~$\p$ on $\{1,2,\ldots,n\}$; we shall suppress
the index~$n$ where possible. The number of weak exceedances is defined to 
be the sum $\sn I_i$, where $I_i := \bone_{\{\p(i) \ge i\}}$.
The process $S_0(t) := \sint I_i$ is thus
of the kind studied in the introduction, with
$a_0(i,j) := \bone_{\{i\le j\}}$. Simple calculations show that
$\ex I_i \Eq \ba_0(i,+) \Eq (n-i+1)/n$, and thus 
\begin{align}
    a(i,j)&\Eq a\un(i,j)\Eq\bone_{\{i\le j\}}-1+(i-1)/n,
        \\
    \ex S_0(k/n)& \Eq \frac{k(2n-k+1)}{2n}.\label{e2}
\intertext{Hence, as $n\to\infty$,}
\label{S0-mean}
   \ex S_0(t) &\Eq nt(1-t/2) + O(1).
\end{align}

Further, although we will not need it,
for $i<j$,
\eqs
   &&\ex\{I_i \giv I_j=1\} \Eq \frac{n-i}{n-1},\qquad 
    \ex\{I_iI_j\} \Eq \frac{(n-i)(n-j+1)}{(n-1)n}, \non\\
\ens
which makes it possible to calculate variances and covariances exactly.
Higher moments can be computed exactly, too.

We now turn to the approximation of
$S(t) := S_0(t) - \ex S_0(t)$. 
We first note that 
$$
   |a(i,j) - \a(i/n,j/n)| \Le n^{-1},
$$
where $\a(t,u) := \bone_{\{t\le u\}} - 1 + t$, so that
$|\a_n(t,u) - \a(t,u)| \le 2n^{-1}$ for $|t-u| > n^{-1}$, and
that $|\a_n(t,u) - \a(t,u)| \le 1$ for all $t,u \in I$. 
Thus $\a_n \to \a$ in~$L_2$, with
\[
   \|\a\|_2^2 \Eq 1/6;\qquad  \|\e_n\|_2^2 \Le 18/n;
   \qquad \a^+ \Eq \sqrt6,
\] 	 
and $\|\a_n\|_\infty/ \|\a\|_\infty$ is bounded.
Calculation based on~\Ref{fun-def-a} 
shows also that, for $0 \le t \le u \le 1$,
\eqs
  f(t) &=& 6\int_0^t x(1-x)\,dx \Eq 3t^2-2t^3; \\
  g(t,u) &=& 6\int_0^t \int_0^u\{(1-x\vee y)- (1-x)(1-y)\}\,dxdy \\
        &=& 3t^2u - t^3 - \tfrac32t^2u^2,
\ens
and that we can take $\b=2$ in~\Ref{Fernique}.
Hence we can apply Theorem~\ref{a-convergent}, and 
defining $Y_n$ by \Ref{Y-def1} with \Ref{s(a)-def},
conclude that
$Y_n \to Z$ in $D[0,1]$, with convergence rate $O\bigl(n^{-1/2}\slogn\bigr)$ 
as measured by~$M_0$-functionals, where~$Z$ is the Gaussian process
given by~\Ref{Z-def-b}:
\[
   Z(t) \Eq \sqrt6 \intiti \{\bone_{\{v\le w\}} - 1 + v\}\,K(dv,dw).
\]
Note also that
\eq\label{Yn-approx}
   Y_n(t) \Eq \sqrt{6/n}\{S_0(t) - nt(1-t/2)\} + O(n^{-1/2}),
\en
indicating that the approximation can be simplified, as in the following
theorem.

\begin{theorem}\label{Th3}
Let $S_0\un(t) := \sint I_i\un$, where $I_i\un := \bone_{\{\p\un(i) \ge i\}}$
and~$\p\un$ is a uniform random permutation on $\{1,2,\ldots,n\}$.
Write $\m(t) := t(1-t/2)$.
Then 
\[
 \hY_n \Def n\qqw\,\{S_0\un - n\m\} \ \to_d\ \hZ \quad \mbox{in}\quad D[0,1],
\]
where $\hZ$ is a zero mean Gaussian process with covariance function~$\hgs$
given by 
\[
  \hgs(t,u) \Eq \tfrac16  \s(t,u) 
    \Eq \tfrac16(f(t) - g(t,u)) \Eq \tfrac12t^2(1-u+\half u^2) - \tfrac16t^3, 
    \quad 0\le t\le u\le 1.
\]
\end{theorem}

The number of weak exceedances of a permutation is one of a number
of statistics that can be deduced from the permutation tableaux
introduced by Steingr\'imsson and Williams~(2007).  Such a tableau 
is a Ferrers diagram (possibly with some rows of length 0) 
with elements
from the set $\{0,1\}$ assigned to the cells, with the following
restrictions:
\begin{enumerate}
\item Each column of the rectangle contains at least one~$1$;
\item There is no~$0$ that has a~$1$ above it in the same column
{\it and\/} a~$1$ to its left in the same row.
\end{enumerate}
The length of a tableau is defined to be the sum of the numbers of its 
rows and columns, and the set of possible tableaux of length~$n$ is in 
one-to-one correspondence with the permutations of~$n$ objects.  In 
particular, under the bijection between tableaux and permutations
defined by Steingr\'imsson and Williams~(2007, Lemma~5), the lower right boundary, 
which consists of a sequence of $n$ unit steps down or to the left, has its 
$i$-th step down if $I_i\un = 1$
and to the left if  $I_i\un = 0$.  Hence the Theorem~\ref{Th3} above, together
with~\Ref{S0-mean}, provides information about the asymptotic shape
of the lower right boundary~$\G_n$ of the tableau corresponding to a randomly
chosen permutation. Let the upper left corner of the Ferrers diagram 
represent the origin with the $x$-axis to the right and the $y$-axis
vertically \emph{downward}, so that the lower right boundary runs from
$(n-S_0(1),0)$ to $(0,S_0(1))$: then $\G_n$ consists of the set
$\{(n-S_0(1)-l + S_0(l),S_0(l)),\,0\le l\le n\}$, linearly interpolated.
Hence, $n\qw\G_n$ is approximated within $O(n\qw)$ by the curve 
\begin{equation*}
  \{(\half[1-t^2] + n\qqw\,(\hY_n(t)-\hY_n(1)),\,
		   \half[1-(1-t)^2] + n\qqw\,\hY_n(t)),
\,0\le t\le 1\},
\end{equation*}
where $\hY_n$	is as defined in Theorem~\ref{Th3}.

\begin{corollary}\label{boundary}
As $n\to\infty$, $n\qw\G_n$ can be approximated in distribution by 
\eqs  
&& \bigl\{(\half[1-t^2] + n\qqw\,(\hZ_n(t)-\hZ_n(1)),\,
		   \half[1-(1-t)^2] +n\qqw\,\hZ_n(t)),\\
&&\hspace{1in}\,0\le t\le 1\bigr\},
\ens
with an error $o(n\qqw)$.
\end{corollary}

In particular, as can also be seen more directly,
$n\qw\G_n$ converges in probability to the deterministic curve
\begin{equation*}
 \bigset{(\half[1-t^2], \half[1-(1-t)^2],\,0\le t\le 1}
=\bigset{(x,y)\in[0,\infty)^2:x+y=\tfrac34-(x-y)^2},
\end{equation*}
an arc of a parabola.

Another statistic of interest is the area~$A_n$ of such a tableau,
which is given by the formula $A_n := \sn I_i \sum_{j=i+1}^n (1-I_j)$,
again because of the bijection above.  Direct computation yields
the expression 
\eqs
   A_n &=& \sn S_0(i/n) - \half S_0^2(1) - \half S_0(1) \\
   &=& \sn \{i(1-i/2n) + \sqrt{n}\,\hY_n(i/n)\}
	   - \half\{(n/2) + \sqrt{n}\,\hY_n(1)\}^2 \\
   &&\qquad\mbox{}- \half\{(n/2) + \sqrt{n}\,\hY_n(1)\}\\
	 &=& \frac{5n^2-2}{24} + n^{3/2}\Bigl\{n^{-1}\sn \hY_n(i/n) - \half\hY_n(1)\Bigr\}\\
   &&\qquad\mbox{}
        -  \half\{\sqrt{n}\,\hY_n(1) + n\hY_n(1)^2\}.
\ens			
This leads to the following limiting approximation.

\begin{corollary}\label{area}
As $n\to\infty$, 
\[
    n^{-3/2}\left(A_n - \frac{5n^2}{24}\right)\ \to_d\ \nn(0,\tfrac1{144}).
\]
\end{corollary}

\proof
By the continuous mapping theorem and Slutsky's lemma, it is immediate 
from Theorem~\ref{Th3} that
\[
   n^{-3/2}\Bigpar{A_n - \frac{5n^2}{24}}\ 
	        \to_d\ \int_0^1 \hZ(t)\,dt - \half \hZ(1).
\]
Now the random variable $\{\int_0^1 \hZ(t)\,dt - \half \hZ(1)\}$ has
mean zero and variance
\[
   \int_0^1 \int_0^1 \hgs(t,u)\,du\,dt - \int_0^1 \hgs(t,1)\,dt + \quarter \hgs(1,1),
\]
with~$\hgs$ as in Theorem~\ref{Th3}, and this gives the value $1/144$.
The corollary follows.
\ep

Note also that the number of rows in the permutation tableau 
$R_n=S_0(1)$;
hence Theorem \ref{Th3} implies also, using $\hgs(1,1)=1/12$,
\[
    n^{-1/2}\left(R_n - \tfrac12 n\right)\ \to_d\ \nn(0,\tfrac1{12}).
\]
This, however, does not require the functional limit theorem; it
follows by the arguments above from 
Hoeffding's (1951) combinatorial central limit theorem, and
it can also be shown in other ways, see 
Hitczenko and Janson (2009+).

\end{document}